\documentclass[12pt]{article}
\usepackage{amsfonts}
\newtheorem{lemma}{Lemma}
\newtheorem{theorem}{Theorem}
\newtheorem{proposition}{Proposition}
\newtheorem{Definition}{Definition}
\newenvironment{proof}{{\bf Proof}.}{\hfill $\Box$}
\newtheorem{remark}{Remark}

\begin{document}

\title{Polynomial extensions of the Weyl $C^*$-algebra}
\author{Luigi Accardi,\\\vspace{-2mm}\scriptsize Volterra Center, University of Roma Tor Vergata\\\vspace{-2mm}\scriptsize Via Columbia 2, 00133 Roma, Italy\\\vspace{-2mm}\scriptsize e-mail: accardi@volterra.uniroma2.it\bigskip\\
Ameur Dhahri\\\vspace{-2mm}\scriptsize Department of Mathematics, Chungbuk National University\\\vspace{-2mm} \scriptsize 1 Chungdae-ro, Seowon-gu, Cheongju, Chungbuk 362-763, Korea\\\vspace{-2mm}\scriptsize e-mail: ameur@chungbuk.ac.kr}
\date{}
\maketitle
\begin{abstract}
We introduce higher order (polynomial) extensions
of the unique (up to isomorphisms) non trivial 
central extension of the Heisenberg algebra, which can be concretely 
realized as sub--Lie algebras of the polynomial algebra generated by 
the creation and annihilation operators in the Schroedinger representation.
The simplest non--trivial of these extensions (the quadratic one) 
is isomorphic to the Galilei algebra, widely studied in quantum physics.\\
By exponentiation of this representation we construct the corresponding 
polynomial analogue of the Weyl $C^*$--algebra and compute the 
polynomial Weyl relations.
From this we deduce the explicit form of the composition law 
of the associated nonlinear extensions of the $1$--dimensional 
Heisenberg group. 
The above results are used to calculate a simple explicit form of the 
vacuum characteristic functions of the non--linear field operators of 
the Galilei algebra, as well as of their moments. The corresponding measures
turn out to be an interpolation family between Gaussian and Meixner, in 
particular Gamma.
\end{abstract}

\section{Introduction}

The program to develop a theory of renormalized higher ($\geq 2$) powers 
of white noise \cite{AcBouk09sg} has led to investigate the 
connections between central extensions and renormalization.
This has, in its turn, led to the discovery of the unique (up to 
isomorphisms) non trivial central extension, denoted $\hbox{Ce--}heis (1)$, 
of the $1$--dimensional Heisenberg algebra $heis(1)$ and of its 
identification with the Galilei algebra (cf \cite{1}, \cite{4}).
In the present paper we will work in the $1$--dimensional Schroedinger 
representation, where $heis (1)$ can be isomorphically realized as the 
real linear span of $\{1,p,q\}$, ($1$ central element) and the operators
$$
(q\varphi)(x):=x\varphi(x);
\quad\quad p\varphi(x):=\quad {1\over i}{\partial \over \partial x}\varphi(x)
\quad ; \quad \varphi \in L^{2} (\mathbb R)
$$
are defined on suitable domains in $L^{2} (\mathbb R)$
and satisfy the Heisenberg commutation relations
\begin{equation}\label{CCR1}
[q,p] = i:= i\cdot 1
\end{equation}
all other commutators being zero. 
In the same representation $\hbox{Ce--}heis (1)$ can be realized as 
the real linear span of
$$
\{1,p,q,q^2\}
$$
Therefore, for each $N\in \mathbb N$, the complex Lie algebra with generators
$$
\{1,p, q, \dots, q^N \}
$$
in the following denoted $heis (1, 1,N)$, is a natural generalization 
of the Galilei algebra.
The analogy with the usual Heisenberg algebra, i.e 
$heis (1,1,1)=:heis (1)$ naturally suggests the following problems:
\begin{enumerate}
\item[1)] To describe the Lie group generated in the Schroedinger representation by $heis (1,1,N)$.
\item[2)] To find an analogue of the Weyl commutation relations for the group in item 1). This is equivalent to describe the $(1,N)$--polynomial extensions of the Heisenberg group.
\item[3)] To determine the vacuum distribution of the nontrivial elements of $heis (1,1,N)$ 
and their moments.
These elements have the form
$$
up +   \sum_{j=0}^N a_j q^j
$$
with $u\neq 0$ and some $a_j\neq 0$ whith $j\geq 2 $.
\end{enumerate}
In the following we solve problems 1) and 2) for general $N$ 
(see see Theorem \ref{composition law} and Proposition \ref{struc-sigmaN} 
respectively) and problem 3) for $N=2$, while for $N>2$ we reduce the
problem to the calculation of a classical expectation value.
A solution of problem 3) was also found in \cite{1} for the 
Galilei algebra ($N=2$) but the technique used there cannot be applied 
to the case $N>2$.
In section \ref{q-proj-meth} we introduce the {\it $q$--projection-method} 
in order to overcome this difficulty. 
This reduces the problem to the calculation of the expectation value 
of functionals of the form $\exp(P(X)+icX)$ where $P$ is a polynomial 
with real coefficients, $c$ is a real number, and $X$ is a standard 
Gaussian random variable. If $P$ has degree $2$, this gives a different 
proof of result obtained in \cite{1}. The calculation of the moments 
in subsection \ref{3.3} are new even in the quadratic case.
The possibility of a continuous (second quantization) version 
of the polynomial extensions of the Weyl $C^*$-algebra is proved 
in the paper \cite{acc-dha2014}.

\section{The $1$--mode $n$--th degree Heisenberg
$*$--Lie algebra $heis_\mathbb{R}(1,n)$ }

\begin{Definition}\label{df-heis(1,n)}
For $n\in\mathbb N^*$ the \textit{$1$--mode $n$--th degree Heisenberg 
algebra},\\ denoted $heis_\mathbb{R}(1,n)$, is the pair
$$
\{V_{n+2},(L_j)^{n+1}_{j=0}\}
$$
where:
\begin{enumerate}
\item[-] $V_{n+2}$ is a $(n+2)$--dimensional real $*$--Lie algebra;
\item[-] $(L_j)^{n+1}_{j=0}$ is a skew--adjoint linear basis of $V_{n+2}$;
\item[-] the Lie brackets among the generators are given by
\begin{eqnarray*}\label{CR-Heis1n1}
[L_i,L_j]=0\qquad ;\qquad \forall i,j\in\{0,1,\cdots,n\}
\end{eqnarray*}
\begin{eqnarray*}\label{CR-Heis1n2}
[L_{n+1},L_k]=kL_{k-1}
\qquad ;\qquad \forall k\in\{1,\cdots,n\}
 \ , \ L_{-1}:=0
\end{eqnarray*}
\end{enumerate}
\end{Definition}
\begin{remark}
\begin{enumerate}
\item[1)] Multiplying each of the generators $(L_j)^{n+1}_{j=0}$ by 
a strictly positive number, one obtains a new basis $(L_j')^{n+1}_{j=0}$ 
of $V_{n+2}$ satisfying the new commutation relations
\begin{eqnarray*}\label{CR-Heis1n3}
[L'_{n+1},L'_k]=kd_kL'_{k-1}
\qquad ;\qquad \forall k\in\{0,\cdots,n\}
 \ , \ L'_{-1}:=0
\end{eqnarray*}
In this case we speak of \textit{a re--scaled copy of} 
the $1$--mode $n$--th degree Heisenberg algebra.
\item[2)] Denoting $\mathbb{R}_n[X]$ the vector space of polynomials
in one indeterminate with real coefficients and degree less
or equal than $n$, the assignment of the basis $(L_j)^{n+1}_{j=0}$
uniquely defines the parametrization
$$
(u,(a_k)_{k\in\{0,1,\dots,n\}})
\in \mathbb{R}\times\mathbb{R}_n[X] \equiv \mathbb R^{n+2} \ \mapsto
 \ \ell_0 (u,P):= \\
$$
\begin{equation}\label{gen-el-1md-alg}
:= uL_{n+1} + \sum_{k\in\{0,1,\dots,n\}}a_kL_{k}
=:uL_{n+1}+P(L) \in heis_\mathbb{R}(1,n)
\end{equation}
of $heis_\mathbb{R}(1,n)$ by elements of $\mathbb{R}_n[X]$.
When no confusion is possible we will use the identification
\begin{equation}\label{id-el-coord}
\ell_0(u,P) \equiv 
(u,(a_k)_{k\in\{0,1,\dots,n\}})
\in \mathbb R^{n+2}
\end{equation}
\end{enumerate}
\end{remark}

\section{The Schroedinger representation of $heis_\mathbb{R}(1,n)$}\label{Schred-repr}

Let $p,q,1$ be the usual momentum, position and identity
operators acting on the one mode boson Fock space
\begin{eqnarray*}\label{Fock}
\mathcal{H}_1=\Gamma(\mathbb{C})=L^2(\mathbb{R})
\end{eqnarray*}
and satisfying the commutation relations (\ref{CCR1}).
The maximal algebraic domain $\mathcal D_{max}$ (see \cite{[AcAyOu03]}),
consisting of the linear combinations of vectors of the form
$$
q^np^k\psi_z
\qquad; \qquad k,n\in\mathbb{N} \ , \ z\in\mathbb{C}
$$
where $\psi_z$ is the exponential vector associated to
$z\in\mathbb{C}$, is a dense subspace of $\Gamma(\mathbb{C})$
invariant under the action of $p$ and $q$ hence of all the
polynomials in the two non commuting variables $p$ and $q$.
In particular, for each $n\in\mathbb N$, the complex linear span of
the set $\{i1,ip,iq,\dots ,iq^n\}$, denoted
$heis_\mathbb{R}(F,1,n)$, leaves invariant
the maximal algebraic domain $\mathcal D_{max}$. Hence the commutators
of elements of this space are well defined on this domain and 
one easily verifies that they define a structure of $*$--Lie algebra
on $heis_\mathbb{R}(F,1,n)$. Since $\mathbb R_n[X]$ has dimension $n+1$, 
this real Lie algebra has dimension $n+2$ on $\mathbb R$.
\begin{Definition}\label{defi1}
For $n\in\mathbb N$, the real Lie algebra 
$$
heis_\mathbb{R}(F,1,n) :=
\big\{iup+iP(q)\quad :\quad u\in\mathbb R,P\in\mathbb R_n[X]\big\}$$
is called {\it the $(1,n)$--polynomial extensions of the $1$--dimensional 
Heisenberg algebra}. 
\end{Definition}

\begin{lemma}  
In the above notations the map  
\begin{eqnarray*}\label{sk-ad-Heis1n-id-gen}
L_{n+1}\mapsto ip \quad , \quad L_0\mapsto i1 \quad , \quad
L_k\mapsto iq^k \quad ; \quad k\in \{1, \dots  , n\}
\end{eqnarray*}
admits a unique linear extension from $heis_\mathbb{R}(1,n)$ onto $heis_\mathbb{R}(F,1,n)$
which is a $*$--Lie algebra isomorphism
called \textit{the Schroedinger representation} of the
$n$--th degree Heisenberg algebra $heis_\mathbb{R}(1,n)$.
\end{lemma}
\begin{proof} The linear space isomorphism property follows from the linear
independence of the set $\{1,p,q,\dots,q^n\}$.
The $*$--Lie algebra isomorphism property follows from direct computation.
\end{proof}

In the following, we will use the Shroedinger representation of the $1$--mode $n$--th degree Heisenberg $*$--Lie algebra and we will identify
$$heis_\mathbb{R}(1,n)\equiv heis_\mathbb{R}(F,1,n)$$

\section{Polynomial extensions of the 1-mode Weyl $C^*$--algebra}

Recall that the position and the momentum operators are defined by
$$
(q\varphi)(x):=x\varphi(x);
\quad\quad p\varphi(x):=\quad {1\over i}{\partial \over \partial x}\varphi(x)
$$
on suitable domains in $L^{2} (\mathbb R)$.  Let $\mathcal{C}_c^\infty(\mathbb{R})$ be the vector space of the $\mathcal{C}^\infty$--functions $f:\mathbb{R}\rightarrow\mathbb{C}$ with compact supports.

For the following lemma we refer the reader to \cite{K}.
\begin{lemma} \label{K} A symmetric operator $A$ defined on a total set of an Hilbert space $\mathcal{H}$ is essentially self-adjoint if and only if $Im(A+i)$ is a total set in $\mathcal{H}$. 
\end{lemma}
In our case, for any $u,w\in\mathbb{R}$ and $P\in\mathbb{R}_n[X]$, one has
$$(wp+uP(q)+i)\big(\mathcal{C}_c^\infty(\mathbb{R})\big)=\mathcal{C}_c^\infty(\mathbb{R})$$
which is a total set in $L^2(\mathbb{R})$. This proves that, for any $u,w\in\mathbb{R}$ and $P\in\mathbb{R}_n[X]$, the operator $(wp+uP(q))$ is essentially self-adjoint. Hence, the operator
$$e^{i(wp+uP(q))}$$
is unitary.

Now, we will use the known fact that, if $\mathcal{L}$ is a Lie algebra 
and $a,b\in\mathcal{L}$ are such that
$$
[a,[a,b]]=0
$$
then, for all $w,z\in\mathbb{R}$
\begin{equation}\label{*1}
e^{za}e^{wb}e^{-za}=e^{w(b+z[a,b])}
\end{equation}
\begin{proposition}\label{prop2}
For all $u,w\in \mathbb{R}$ ($w\neq0$) and all $P\in\mathbb{R}_n[X]$, we have
\begin{equation}\label{*2} 
e^{iwp+iuP'(q)}=e^{iu\frac{P(q+w)-P(q)}{w}}e^{iwp}
\end{equation}
\end{proposition}
\begin{proof}
Let $P$ a polynomial in one determinate and $w,z\in\mathbb{R}$ such that $w\neq0$. Then from the identity 
$$
[P(q),p]=iP'(q)
$$
It follows that
$$
[P(q),[P(q),p]]=0
$$
Hence (\ref{*1}) implies that
\begin{eqnarray}\label{pq}
e^{izP(q)}e^{iwp}e^{-izP(q)}=e^{iw(p-zP'(q))}
\end{eqnarray}
But, one has also
\begin{eqnarray}\label{pp}
e^{iwp}e^{-izP(q)}=(e^{iwp}e^{-izP(q)}e^{-iwp})e^{iwp}=e^{-izP(q+w)}e^{iwp}
\end{eqnarray}
The identities (\ref{pq}) and (\ref{pp}) imply 
\begin{equation}\label{pqp}
e^{iw(p-zP'(q))}= e^{izP(q)}e^{-izP(q+w)}e^{iwp}
\end{equation}
With the change of variable $u:=-wz$, (\ref{pqp}) becomes (\ref{*2}).
\end{proof}

In the following we will use the identification
$$
P'(X)=\sum_{k=0}^na_kX^k\equiv(a_0,a_1,\cdots,a_n)\in\mathbb R^{n+1}
$$
Now, for all $w\in\mathbb{R}$, we define the linear map $T_w$ by 
$T_0:=id$ if $w=0$ and if $w\neq0$ 
\begin{equation}\label{3.1}
T_w:P'\in\mathbb R_n[X]\rightarrow 
\frac{1}{w}\int_0^wP'(.+y)dy\in \mathbb R_n[X]
\end{equation}
so that 
\begin{equation}\label{3.2} 
T_wP'={{P(\cdot+w)-P}\over w}
\end{equation}
where $P$ is any primitive of $P'$.
\subsection{Properties of the maps $T_w$}
\begin{lemma}\label{T_w}
For each $w\in\mathbb R\setminus\{0\}$, the matrix of the linear map 
$T_w$ with respect to the mononial real basis $\{1,X,\dots,X^n\}$ 
of $\mathbb R_n[X]$, still denoted by $T_w$, has the form
\begin{equation}\label{Tw=1+T(X,w)} 
T_w= 1+T_{X,w}
\end{equation}
where $T_{X,w}$ is a strictly lower triangular nilpotent matrix 
whose coefficients are given by:
\begin{equation}\label{tmn}
t_{m_1m_2}(w)
=\chi_{\{m_1<m_2\}}\frac{m_2!}{(m_2+1-m_1)!m_1!}w^{m_2-m_1}
\quad ;\quad 0\leq m_1,m_2\leq n
\end{equation} 
$$
\chi_{\{m_1<m_2\}}=\left\{
\begin{array}{cc}
0,& \mbox{ if }\,m_1\geq m_2\\
1,& \mbox{ if }\, m_1<m_2
\end{array}
\right.
$$
In particular $T_w$ is invertible with inverse given by 
\begin{equation}\label{inv-Tw}
T_w^{-1}=(I+T_{X,w})^{-1}=\sum_{k=0}^{n}(-1)^k(T_{X,w})^k
\end{equation}
\end{lemma}
\begin{proof} Clearly $T_w$ is linear and 
\begin{equation}\label{Tw1=1}
T_w1=1
\end{equation}
Let $k\in\{0,1,\cdots,n\}$. If we take 
$$
P'(X):=X^k
$$
then, in (\ref{3.2}) we can choose
$$
P(X)={{X^{k+1}}\over{k+1}}
$$
so that
$$
{{P(X+w)-P(X)}\over w}={1\over{w(k+1)}}[(X+w)^{k+1}-X^{k+1}]
$$
Since
\begin{eqnarray*}
(X+w)^{k+1}&=&\sum^{k+1}_{h=0}\pmatrix{k+1\cr h\cr}X^h w^{k+1-h}\\
&=&X^{k+1}+\sum^k_{h=0}\pmatrix{k+1\cr h\cr}X^h w^{k+1-h}
\end{eqnarray*}
it follows that
$$
(X+w)^{k+1}-X^{k+1}=\sum^k_{h=0}\pmatrix{k+1\cr h\cr}X^hw^{k+1-h}
$$
Hence for all $k\geq1$
\begin{eqnarray}\label{explicit}
T_w X^k&=&\frac{1}{w(k+1)}\sum_{h=0}^k\pmatrix{k+1\cr h\cr}w^{k+1-h}X^h\nonumber\\
&=&\sum_{h=0}^k\frac{k!}{(k+1-h)!h!}w^{k-h}X^h\nonumber\\
&=&\sum_{h=0}^{k-1}\frac{k!}{(k+1-h)!h!}w^{k-h}X^h+X^k
\end{eqnarray}
Therefore the matrix of $T_w$ in the basis $\{1,X,\dots,X^n\}$ of 
$\mathbb{R}_n[X]$ has the form (\ref{Tw=1+T(X,w)}) with the entries 
$(t_{m_1m_2}(w))$ of $T_{X,w}$ given by (\ref{tmn}).
Thus, in the same basis, $T_{X,w}$ is a strictly lower triangular 
nilpotent matrix and it is known that this implies that 
$T_w$ is invertible with inverse given by (\ref{inv-Tw}).
\end{proof} \\

Let $S$ denote the action of $\mathbb{R}$ on $\mathbb{R}_n[X]$ by translations
$$(S_uP')(X):=P'(X+u),\;X,\,u\in\mathbb{R}$$
Then 
$$
T_uP'={1\over u}(S_u-id)P
$$
where $P$ is a primitive of $P'$. 
The shift $S_u$ is a linear homogeneous map of $\mathbb R_n[X]$ 
into itself. Moreover, with $P$ and $P'$ as above, one has:
$$
S_u(P(X))=
\sum^n_{k=0}a_k(X+u)^k
=\sum^n_{k=0}\sum^k_{h=0}a_k\pmatrix{k\cr h\cr}X^h u^{k-h}=
$$
$$
=\sum^n_{h=0}\Big(\sum^n_{k=h}a_k\pmatrix{k\cr h\cr}u^{k-h}\Big)X^h
=\sum^n_{h=0}(S_uP')_hX^h
$$
where 
$$
(S_uP')_h:=\sum^n_{k=h}a_k\pmatrix{k\cr h\cr}u^{k-h}
$$
Thus the matrix of $S_u$ in the monomial basis $\{1,X,\dots,X^n\}$ 
of $\mathbb R_n[X]$,  is
\begin{eqnarray}\label{translation}
(S_u)_{hk}=\chi_{\{h\leq k\}}\pmatrix{k\cr h\cr}u^{k-h}
\end{eqnarray}
Clearly $S_u$ is invertible with inverse $S_{-u}$ and
$$
S_uS_v=S_{u+v}
$$
Denote $\partial$ the derivation operator on $\mathbb{R}_n[X]$ defined by
\begin{equation}\label{d-par}
\partial:\sum_{k=0}^na_kX^k\rightarrow\sum_{k=1}^nka_kX^{k-1}
\end{equation}
Then clearly 
\begin{equation}\label{ds=sd}
\partial S_u=S_u\partial\qquad ;\qquad \forall u\in\mathbb{R}
\end{equation}
Let $\mathbb{R}_n[X]_0$ the subalgebra of $\mathbb{R}_n[X]$ of all polynomials vanishing in zero
$$
\mathbb{R}_n[X]_0:=\{P'\in \mathbb{R}_n[X]:\;P'(0)=0\}
$$
Then
$$
\mathbb{R}_n[X]=\mathbb{R}_n[X]_0\oplus\mathbb{R}1
$$
and, thanks to (\ref{Tw1=1}), $T_w$ is uniquely determined by 
its restriction to $\mathbb{R}_n[X]_0$. Denote
\begin{eqnarray*}
T_w^0&:=&T_w\big|_{\mathbb{R}_n[X]_0}\\
\partial_0&:=&\partial\big|_{\mathbb{R}_n[X]_0}
\end{eqnarray*}
Then (\ref{d-par}) implies that $\partial_0$ is invertible and from (\ref{ds=sd}) one has
\begin{equation}\label{IS=SI}
\partial_0^{-1}S_u=S_u\partial_0^{-1},\;\forall u\in\mathbb{R}
\end{equation}
\begin{lemma}\label{lem2}
The following identities hold:
\begin{enumerate}
\item[1)] $T_wS_u=S_uT_w,\;\forall u,w\in\mathbb{R}$
\item[2)] $T_wT_u=T_uT_w,\;\forall u,w\in\mathbb{R}$
\item[3)] $T_uS_v=(1+\frac{v}{u})T_{u+v}-\frac{v}{u}T_v
\quad ;\quad \forall u,v\in\mathbb{R}$ such that $u\neq0$
\item[4)] $T_uS_{-u}=T_{-u},\;\forall u\in\mathbb{R}$
\item[5)] $T_uT_v=(\frac{1}{v}+\frac{1}{u})T_{u+v}
-\frac{1}{u}T_v-\frac{1}{v}T_u, \;\forall u,v\in\mathbb{R}\setminus\{0\}$
\end{enumerate}
\end{lemma}
\begin{proof} In the above notations
$$T_w^0=\frac{1}{u}(S_u-id)\partial_0^{-1}$$
and (\ref{IS=SI}) implies that
$$T_w^0S_u=S_uT_w^0,\;\forall u,w\in\mathbb{R}$$
Since the constant function $1$ is fixed both for $S_u$ and for $T_w$, 
it follows that
$$
T_wS_u=S_uT_w,\;\forall u,w\in\mathbb{R}
$$
and therefore also
$$
T_wT_u=T_uT_w,\;\forall u,w\in\mathbb{R}
$$
which implies that assertions 1) and 2) are satisfied.\\
Now notice that translations commute with the operation of taking 
primitives, i.e. if $P$ is a primitive of $P'$, then 
$\forall v\in\mathbb{R},\;S_vP$ is a primitive of $S_vP'$, that is
\begin{eqnarray*}
T_uS_vP'&=&\frac{1}{u}(S_u-id)S_vP\\
&=&\frac{1}{u}(S_{u+v}-S_v)P\\
&=&\frac{1}{u}(S_{u+v}-id)P+\frac{1}{u}(id-S_v)P\\
&=&\Big(\frac{u+v}{u}T_{u+v}-\frac{v}{u}T_v\Big)P'
\end{eqnarray*}
This proves 3). In particular, choosing $v=-u$, one finds 
$T_uS_{-u}=T_{-u}$ and 4) holds. Finally, from  the identity
$$
T_uT_v=\frac{1}{v}T_uS_v-\frac{1}{v}T_u
$$
and from 3), the assertion 5) follows.
\end{proof}

\subsection{The Polynomial extensions Heisenberg group}

Recall that the $1$--dimensional Heisenberg group, denoted $Heis(1,1)$,
is the nilpotent Lie group whose underlying manifold is 
$\mathbb R\times \mathbb R^2$ and whose group law is given by
$$
(t,z)\circ (t',z') = (t+t'-2\sigma (z, z') ,z+z'); \quad
t,t'\in\mathbb R \ ; \ z:=(\alpha,\beta),z':=(\alpha',\beta')\in\mathbb R^2
$$
where $\sigma ( \ \cdot \ , \ \cdot \ ) $ the symplectic form on 
$\mathbb R^2$ given by
\begin{equation}\label{sympl-form}
\sigma (z, z') := \alpha\beta'-\alpha'\beta\qquad ; \qquad 
z:=(\alpha,\beta),z':=(\alpha',\beta')\in\mathbb R^2
\end{equation}
In the Schroedinger representation of the real Lie algebra 
$heis_\mathbb{R}(1,1)$ the {\it Weyl operators} are defined by
$$
W(z) :=e^{ i \sqrt 2 (\beta q-\alpha p) }\qquad ; \qquad z:=(\alpha,\beta)\in\mathbb R^2
$$
and the {\it extended Weyl operators} by
\begin{equation}\label{un-rep1}
W(t,z) := e^{ i (\sqrt 2 (\beta q-\alpha p)  + t1) } = W(z)e^{ it }
\qquad ; \qquad (t,z)\in\mathbb R\times \mathbb R^2
\end{equation}
From the Weyl relations
\begin{equation}\label{cl-Weyl-rel}
W(u)W(z)  = e^{- i\sigma (u, z) } W( z  +  u)  
\end{equation}
one then deduces that
\begin{equation}\label{un-rep1}
W(t,z)W(t',z') = W((t,z)\circ (t',z'))
\end{equation}
i.e. that the map (\ref{un-rep1}) gives a unitary representation of the group $Heis(1,1)$.\\
The generators of the centrally extended Weyl operators have the form
\begin{equation}\label{upP'q11}
up+a_1q + a_01 := up+P'(q) 
\end{equation}
where $u\in\mathbb R$ and $P'$ is a polynomial in one indeterminate 
with real coefficients of degree at most $1$. Thus they are exactly
the elements of the one dimensional Heisenberg algebra $heis_\mathbb{R}(1,1)$.


With these notations the $(1,n)$--polynomial extensions 
of the centrally extended Weyl operators are the operators of the form
\begin{equation}\label{df-1Npol-cent-ext-W-op}
W(u,P'):=e^{i(up+P'(q)) }
\quad ; \ u\in\mathbb R \quad ; \ P'\in \mathbb{R}_n[X]
\end{equation}
By analogy with the $1$--dimensional Heisenberg group, we expect that
the pairs $(u,P')\in \mathbb R \times \mathbb{R}_n[X]$ form a group for 
an appropriately defined composition law. The following theorem shows that 
this is indeed the case.
\begin{theorem}\label{composition law}  
For any $(u,P'),(v,Q')\in \mathbb R \times \mathbb{R}_n[X]$ one has
\begin{equation}\label{1Npol-centr-ext-Weyl-rel}
W(u,P')W(v,Q')= W((u,P')\circ (v,Q'))
\end{equation}
where, in the notaions of Lemma (\ref{lem2}):
\begin{equation}\label{1Npol-comp-law1} 
(u,P')\circ(v,Q'):= (u+v,T^{-1}_{u+v}(T_uP'+T_vS_uQ'))
\end{equation}
\end{theorem}
\begin{proof}
From Proposition \ref{prop2}, one has
$$
W(u,P')= e^{iup+iP'(q)}=e^{iT_uP'(q)}e^{iup}
$$
Therefore, using the identity (\ref{*2}) and recalling the notation 
(\ref{3.2}), one deduces that
\begin{eqnarray*}
W(u,P')W(v,Q')&= &e^{iup+iP'(q)}e^{ivp+iQ'(q)}\\
&=&e^{iT_uP'(q)}e^{iup}e^{iT_vQ'(q)}e^{ivp}\\
&=&e^{iT_uP'(q)}e^{iT_vQ'(q+u)}e^{i(u+v)p}\\
&=&e^{i[T_uP'(q)+T_vQ'(q+u)]}e^{i(u+v)p}\\
&=&e^{i(u+v)p+iT^{-1}_{u+v}(T_uP'(q)+T_vS_uQ'(q))}\\
&=&W\left(u+v, T^{-1}_{u+v}(T_uP'+T_vS_uQ')\right)
\end{eqnarray*}
This proves (\ref{1Npol-centr-ext-Weyl-rel}) and (\ref{1Npol-comp-law1}).
\end{proof}

\begin{remark}
The associativity of the group law (\ref{1Npol-comp-law1}) follows 
from the associativity of the multiplication of bounded operators 
on Hilbert spaces.\\ It can be directly verified using Lemma \ref{lem2}.
\end{remark}

Our next goal is to determine the $(1,n)$--polynomial extensions of the 
Weyl commutation relations. To this goal define the ideal
\begin{equation}\label{df-RNx0}
 \mathbb{R}_n[X]_0
:=\{P_0\in\mathbb{R}_n[X] \ : \ P_0(0)=0\}
\end{equation}
Thus the projection map from $\mathbb{R}_n[X]$ onto $\mathbb{R}_n[X]_0$ 
is given by
\begin{equation}\label{df-proj-RNx0}
P'(X)\in \mathbb{R}_n[X]\mapsto (P')_0(X):= P'(X) - P'(0)\in\mathbb{R}_n[X]_0 
\end{equation}
and, in the notation (\ref{df-proj-RNx0}), the $(1,n)$--polynomial 
extensions of the centrally extended Weyl operators 
(\ref{df-1Npol-cent-ext-W-op}) take the form
\begin{equation}\label{df-1Npol-cent-ext-W-op2}
W(u,P')=e^{i(up+P'(q)) }=e^{i(up+(P')_0(q)) }e^{iP'(0) }
\quad ; \ u\in\mathbb R \quad ; \ P'\in \mathbb{R}_n[X]
\end{equation}
The analogy between (\ref{df-1Npol-cent-ext-W-op2}) and (\ref{un-rep1})
naturally suggests the following definition.
\begin{Definition}
The unitary operators 
\begin{equation}\label{df-1Npol-ext-W-op}
W(u,P'_0):=e^{i(up+P'_0(q))}
\quad ;\quad u\in\mathbb R  \quad ;\quad  P'_0\in \mathbb{R}_n[X]_0
\end{equation}
will be called {\it $(1,n)$--polynomially extended Weyl operators}.
\end{Definition}
From  (\ref{1Npol-centr-ext-Weyl-rel}) and (\ref{1Npol-comp-law1}) 
we see that, if $P'_0,Q'_0\in \mathbb{R}_n[X]_0$, then in the 
notation (\ref{df-proj-RNx0}) one has
\begin{eqnarray}\label{27a}
&&W(u,P'_0)W(v,Q'_0)  \nonumber\\
&&=  W\left(u+v, T^{-1}_{u+v}(T_uP_0'+T_vS_uQ_0'\right)\\ 
                    &&= W\left(u+v, (T^{-1}_{u+v}(T_uP_0'+T_vS_uQ_0'))_0\right) 
e^{i(T^{-1}_{u+v}(T_uP_0'+T_vS_uQ_0'))(0)}\nonumber
\end{eqnarray}
By analogy with the usual Weyl relations we introduce the notation
\begin{equation}\label{df-sigmaN}
\sigma ((u,P_0'),(v,Q_0')) :=   
T^{-1}_{u+v}(T_uP_0'+T_vS_uQ_0'))(0)
\end{equation}
i.e. $\sigma ((u,P_0'),(v,Q_0'))$ is the degree zero coefficient of 
the polynomial\\ $T_{u+v}^{-1}(T_uP'_1+T_vS_uP'_2)\in \mathbb{R}_n[X]$.
Notice that the map 
$$
((u,P_0'),(v,Q_0'))\equiv ((u,v),(Q_0',P_0'))\mapsto 
\sigma ((u,P_0'),(v,Q_0')) 
$$
is linear in the pair $(Q_0',P_0')$ but polynomial in the pair $(u,v)$.
This is an effect of the duality between $p$, which appears at the 
first power, and $q$, which appears in polynomial expressions.\\
In order to prove the $(1,n)$--polynomial analogue of the classical 
Weyl relations (\ref{cl-Weyl-rel}) one has to compute the scalar factor 
in (\ref{27a}).
We will compute more generally all the coefficients of
$T_{u+v}^{-1}(T_uP'_1+T_vS_uP'_2)$. In view of (\ref{inv-Tw}) this leads to 
compute the powers of the matrices $T_{X,w}$ given by (\ref{tmn}).
\begin{lemma}\label{LeM}
Let $k\in\{2,\dots,n\}$ and $w\in\mathbb{R}\setminus \{0\}$. 
Then, the matrix of $(T_{X,w})^k$, in the monomial basis 
$\{1,X,\dots,X^n\}$ of $\mathbb{R}_n[X]$, is given by \\
$\big(T^{[k]}_{ij}(w)\big)_{0\leq i,j\leq n}$, where
\begin{equation}\label{k-power}
T^{[k]}_{ij}(w)=t_{ij}(w)C_{i,j}^{[k]}
\end{equation}
the coefficients $t_{ij}(w)$ are given by (\ref{tmn}) and the $C_{i,j}^{[k]}$ are inductively defined by
\begin{equation}\label{28a}
C_{i,j}^{[k]}=\sum_{h=0}^nC_{i,j}^{(h)}C_{i,h}^{[k-1]}
\end{equation}
where
\begin{equation}\label{df-Ci,j(h)}
C_{i,j}^{(h)}=\frac{(j+1-i)!}{(h+1-i)!(j+1-h)!}\chi_{\{i<h<j\}}(h)
\end{equation}
and $C_{i,j}^{[1]}=\chi_{\{i<j\}}$ for all $i,j=0,\dots, n$.
\end{lemma}
\begin{proof}
We prove the lemma by induction. For $k=2$, using (\ref{tmn}) one has
\begin{eqnarray}\label{tihthj}
t_{ih}(w)t_{hj}(w)&=&\chi_{\{i<h\}}\chi_{\{h<j\}}\frac{1}{h+1}\pmatrix{h+1\cr i\cr}w^{h-i}\frac{1}{j+1}\pmatrix{j+1\cr h\cr}w^{j-h}\nonumber\\
&=&\Big[\chi_{\{i<h<j\}}\frac{j!}{(j+1-i)!i!}w^{j-i}\Big]\frac{(j+1-i)!}{(h+1-i)!(j+1-h)!}\nonumber\\
&=&t_{ij}(w)C^{(h)}_{i,j}
\end{eqnarray}
where $C^{(h)}_{i,j}$ is given by (\ref{df-Ci,j(h)}).
It follows that
$$T_{ij}^{[2]}(w)=t_{ij}(w)\sum_{h=0}^nC_{i,j}^{(h)}=t_{ij}(w)C^{[2]}_{i,j}$$
Suppose by induction that (\ref{k-power}) is true for $2\leq k\leq n-1$. 
Then one has
$$
T_{ij}^{[k+1]}(w)
=\sum_{h=0}^nT_{ih}^{[k]}(w)t_{hj}(w)
=\sum_{h=0}^nT_{ih}^{[k]}(w)t_{hj}(w)
=\sum_{h=0}^nt_{ih}(w)t_{hj}(w)C_{i,h}^{[k]}
$$
From (\ref{tihthj}) it then follows that
$$
T_{ij}^{[k+1]}(w)
=t_{ij}(w)\sum_{h=0}^nC_{i,j}^{(h)}C_{i,h}^{[k]}=t_{ij}(w)C_{i,j}^{[k+1]}
$$ 
and this proves the induction step.
\end{proof}\\
As a consequence of Theorem \ref{composition law} and Lemma \ref{LeM}, 
we are able to explicitly compute the scalar factor in the 
$(1,n)$-polynomial extensions of the Weyl relations (\ref{27a}).
\begin{proposition}\label{struc-sigmaN}
In the notation (\ref{df-sigmaN}), for any $u,v\in\mathbb{R}$ 
and for any $P_1',P'_2\in\mathbb{R}_n[X]_0$ given by
$$
P'_1(X)=\alpha_1 X+\dots+\alpha_n X^n\quad ;\quad 
P'_2(X)=\beta_1 X+\dots+\beta_n X^n
$$ 
one has
\begin{eqnarray*}
\sigma ((u,P_0'),(v,Q_0'))&=&\sum_{j=0}^n\Big[\Big(\frac{1}{j+1}(u+v)^j
\sum_{m=0}^n(-1)^mC_{0,j}^{[m]}\Big)\\
&&\;\;\;\;\;\;\Big(\sum_{j\leq h\leq n}\Big\{t_{jh}(u)\alpha_h+t_{jh}(v)\sum_{h\leq k\leq  n}(S_u)_{hk}\beta_k\Big\}\Big)\Big] 
\end{eqnarray*}
where the coefficients $(S_u)_{hk}$ are given by (\ref{translation}) and with the convention 
$$
C_{0,j}^{[0]}=\delta_{j,0},\;\alpha_0=\beta_0=0,\;t_{jj}(w)=1,\;\mbox{ for all } j=0,\dots,n.$$
\end{proposition}
\begin{proof}
Let $u,v\in\mathbb{R}$ and $P_1',P'_2\in\mathbb{R}_n[X]_0$ be as in the statement. Then, using together Lemmas \ref{T_w}, \ref{LeM} and identity (\ref{translation}), we show that 
\begin{eqnarray*}
T_{u+v}^{-1}(T_uP_1'+T_vS_uP'_2)(X)&=&\!\!\!(T_{u+v}^{-1}(T_uP'_1+T_vS_uP'_2))_0(X)\\
&&\!\!\!\!\!\!\!\!+\sum_{j=0}^n\Big[\Big(\frac{1}{j+1}(u+v)^j\sum_{m=0}^n(-1)^mC_{0,j}^{[m]}\Big)\\
\;\;\;\;\;\;\;\;\;&&\Big(\sum_{j\leq h\leq n}\Big\{t_{jh}(u)\alpha_h+t_{jh}(v)\sum_{h\leq k\leq n}(S_u)_{hk}\beta_k\Big\}\Big)\Big]
\end{eqnarray*}
\end{proof}

\begin{Definition}\label{defi3}
For $n\in\mathbb N$, the real Lie group with manifold 
\begin{equation}\label{24a}
 \mathbb R \times \mathbb R\times\mathbb{R}_n[X]_0 \equiv \mathbb R \times \mathbb R_n[X]
\end{equation}
and with composition law given by 

$$ (u,P')\circ(v,Q'):= (u+v,T^{-1}_{u+v}(T_uP'+T_vS_uQ'))$$
is called {\it the $(1,n)$--polynomial extensions of the $1$--dimensional 
Heisenberg group} and denoted $Heis(1,n)$. 
\end{Definition}
\begin{remark}
The left hand side of (\ref{24a}) emphasiyes that the coordinates \\$(t,u,a_1,\dots,a_n)\in\mathbb{R}^{n+2}$ of an element of $Heis(1,n)$ are intuitively interpreted as time $t$, momentum $u$ and coordinates of the first $n$ powers of position $(a_1,\dots,a_n)$.\\ Putting $t=a_0$ is equivalent to the identification $\mathbb{R}\times\mathbb{R}_n[X]_0\equiv \mathbb R_n[X].$
\end{remark}

\subsection{The Heisenberg algebra $\{i1,ip,iq\}$}

For $N=1$ $heis_\mathbb{R}(F,1,n)$ is reduced to the usual Heisenberg 
algebra  $heis_\mathbb{R}(1,1)\equiv\{i1,ip,iq\}$ and the following 
Proposition shows that the non--linear term $\sigma ( \ \cdot \ , \ \cdot \ )$
, given by (\ref{df-sigmaN}), is given by the usual symplectic form 
(\ref{sympl-form}).
This proves that the composition law given in Definition \ref{defi3} is indeed a generalization of the composition law of the Heisenberg group.
\begin{proposition}\label{Weyl}
Let $P'_1(X)=\alpha_1X+\alpha_0,\;P'_2(X)=\beta_1X+\beta_0$ and 
let $u,v\in\mathbb{R}$. Then, we have
$$
(u,P'_1)\circ(v,P'_2)=(u,\alpha_1X+\alpha_0)\circ(v,\beta_1X+\beta_0)=
$$
$$
=\left(u+v, (\alpha_1+\beta_1)X+ (\alpha_0+\beta_0)
+\frac{1}{2}(u\beta_1-v\alpha_1)\right)
$$
\end{proposition}
\begin{proof}
Let $u,\,w\in\mathbb{R}$. Then, one has
\begin{eqnarray*}
T_w=\pmatrix{1 & \frac{1}{2}w \cr
0 & 1\cr} \quad ;\quad 
T_w^{-1}=\pmatrix{1 & -\frac{1}{2}w \cr
0 & 1\cr}
  \quad ;\quad 
S_u=\pmatrix{1 & u \cr 0 & 1\cr}
\end{eqnarray*}
Let $P'_1(X)=\alpha_0+\alpha_1X,\;P'_2(X)=\beta_0+\beta_1X$. 
Then, from Theorem \ref{composition law}, the composition law is given by
$$ 
(u,P'_1)\circ(v,P'_2)=(u+v,Q')
$$
where 
\begin{eqnarray*}
Q'(X)&=&T^{-1}_{u+v}(T_uP'_1+T_vS_uP'_2)(X)\equiv
\end{eqnarray*}
$$
\equiv\pmatrix{1 & -\frac{1}{2}(u+v) \cr 0 & 1\cr}\left(
\pmatrix{1 & \frac{1}{2}u \cr 0 & 1\cr}\pmatrix{\alpha_0 \cr \alpha_1\cr}
+\pmatrix{1 & \frac{1}{2}v \cr 0 & 1\cr}\pmatrix{1 & u \cr 0 & 1\cr}
\pmatrix{\beta_0 \cr \beta_1\cr}\right)=
$$
$$
=
\pmatrix{1 & -\frac{1}{2}v \cr 0 & 1\cr}\pmatrix{\alpha_0 \cr \alpha_1\cr}
+\pmatrix{1 & -\frac{1}{2}u \cr 0 & 1\cr}\pmatrix{1 & u \cr 0 & 1\cr}
\pmatrix{\beta_0 \cr \beta_1\cr}
$$
$$
=
\pmatrix{\alpha_0 -\frac{1}{2}v\alpha_1 \cr \alpha_1\cr}
+\pmatrix{1 & \frac{1}{2}u \cr 0 & 1\cr}
\pmatrix{\beta_0 \cr \beta_1\cr}
$$
$$
=\pmatrix{\alpha_0 -\frac{1}{2}v\alpha_1 \cr \alpha_1\cr}
+\pmatrix{\beta_0 +\frac{1}{2}u\beta_1\cr \beta_1\cr}
=\pmatrix{(\alpha_0+\beta_0 ) +\frac{1}{2}(u\beta_1-v\alpha_1) \cr 
\alpha_1+\beta_1\cr}\equiv
$$
$$
\equiv(\alpha_0+\beta_0+\frac{1}{2}u\beta_1-\frac{1}{2}v\alpha_1)
+(\alpha_1+\beta_1)X
$$
\end{proof}

\subsection{The Galilei algebra $\{i1,ip,iq,iq^2\}$}

In the case of $n=2$, the composition law is given by the following.
\begin{proposition}\label{N=2}
Let $P'_1,\,P'_2\in\mathbb{R}_2[X]$, such that $P'_1(X)=\alpha_0+\alpha_1X+\alpha_2X^2$ and $P'_2(X)=\beta_0+\beta_1X+\beta_2X^2$. Then, for all $u,v\in\mathbb{R}$, we have
$$(u,P'_1)\circ(v,P'_2)=(u+v,Q')$$
where $Q'(X)=\gamma+\beta X+\alpha X^2$ with
\begin{eqnarray*}
\gamma&=&\alpha_0+\beta_0-\frac{1}{2}v\alpha_1+\frac{1}{2}u\beta_1+\frac{1}{6}(v-u)v\alpha_2+\frac{1}{6}(u-v)u\beta_2\\
\beta&=&\alpha_1+\beta_1-v\alpha_2+u\beta_2\\
\alpha&=&\alpha_2+\beta_2
\end{eqnarray*}
\end{proposition}
\begin{proof}
From Lemma \ref{T_w}, the matrices of $T_w$ and $T_w^{-1}$, in the 
monomial basis $\{1,X,X^2\}$ of $\mathbb{R}_2[X]$, are given by
\begin{eqnarray*}
T_w=\left(
\begin{array}{lcc}
1 & \frac{1}{2}w & \frac{1}{3}w^2\\
0& 1& w\\
0& 0 &1
\end{array}
\right),\;\;
T_w^{-1}=\left(
\begin{array}{lcc}
1 & -\frac{1}{2}w& \frac{1}{6}w^2\\
0& 1& -w\\
0& 0 &1
\end{array}
\right)
\end{eqnarray*}
Moreover, identity (\ref{translation}) implies that
$$S_u=\left(
\begin{array}{lcc}
1 & u & u^2\\
0 & 1& 2u\\
0 & 0 &1
\end{array}
\right)
$$
Let $P'_1(X)=\alpha_0+\alpha_1X+\alpha_2X^2$, $P'_2(X)=\beta_0+\beta_1X+\beta_2X^2$ and $u,v\in\mathbb{R}$. Then, using the explicit form of $T_{u+v}^{-1}, \,T_u$ and $S_u$, it is straightforward to show that
$$
T^{-1}_{u+v}(T_uP'_1+T_vS_uP'_2)(X)=\gamma+\beta X+\alpha X^2
$$
where $\alpha,\beta,\gamma$ are given in the above proposition. 
Finally, from Theorem \ref{composition law} we can conclude.
\end{proof}

\section{ Vacuum expectations and the q-projection method }\label{q-proj-meth}

Because of the relation (see Proposition \ref{prop2})
\begin{equation}\label{ccr4}
 e^{i\alpha(\sqrt{2}p)} e^{i\beta(\sqrt{2}q)}
  =  e^{i2\alpha\beta}e^{i\beta(\sqrt{2}q)}e^{i\alpha(\sqrt{2}p)}
\end{equation}
the Weyl algebra coincides with the complex linear span of the products
$e^{i\beta(\sqrt{2}q)}e^{i\alpha(\sqrt{2}p)}$. Therefore a state on the Weyl
algebra is completely determined by the expectation values of these products.
In particular the Fock state $\Phi$ is characterized by the property that the
vacuum distribution of the position operator is the standard Gaussian
$$
\sqrt{2}q \ \sim \ \mathcal{N}(0,1)
$$ 
Note that
\begin{eqnarray*}\label{decomp4}
e^{i\alpha\sqrt{2}p}\Phi&=&e^{\alpha (b^+-ib)}\Phi\\
&=&e^{-\frac{\alpha^2}{2}}e^{\alpha b^+}e^{\alpha b}\Phi\\
&=&e^{-\alpha^2}e^{\sqrt{2}\alpha\frac{(b^++b)}{\sqrt{2}}}\Phi\\
&=&e^{-\alpha^2}e^{\sqrt{2}\alpha q}\Phi,
\end{eqnarray*}
Moreover, from Proposition \ref{prop2}, one has
$$e^{i\alpha\sqrt{2}p+\alpha\sqrt{2} q}=e^{-\alpha^2}e^{-\alpha\sqrt{2} q}e^{i\alpha \sqrt{2}p}$$

The q-projection method consists in using (\ref{decomp4}) and (\ref{pp}) 
to reduce the problem to compute vacuum expectation values of products of 
the form $e^{-izP(q)}e^{-iwp}$ to the calculation of a single Gaussian integral.\\
In the following sub--sections we illustrate this method starting from 
the simplest examples.

\subsection{Vacuum characteristic functions of observables in $Heis(1,1)$}

The $q$--projection method, applied to $Heis(1,1)$, gives the standard 
result for the spectral measure of the Weyl operators. 
In fact from (\ref{decomp4}) and the CCR it follows that
$$
e^{i(\alpha(\sqrt{2}p)+\beta(\sqrt{2}q))}\Phi
=e^{i\alpha\beta}e^{-\alpha^2}e^{(i\beta-\alpha)(\sqrt{2} q)}\Phi
$$
from which one obtains
$$
\langle \Phi,e^{i(\alpha(\sqrt{2}p)+\beta(\sqrt{2}q))}\Phi\rangle
=e^{i\alpha\beta}e^{-\alpha^2}\langle\Phi,e^{(i\beta-\alpha)(\sqrt{2} q)}\Phi\rangle
=\frac{e^{-\alpha^2}e^{i\alpha\beta}}{\sqrt{2\pi}}
\int_\mathbb{R}e^{i\beta x-\alpha x}e^{\frac{-x^2}{2}}dx
$$
$$
=\frac{e^{-\frac{\alpha^2}{2}}e^{i\alpha\beta}}{\sqrt{2\pi}}
\int_\mathbb{R}e^{i\beta x}e^{\frac{-(x+\alpha)^2}{2}}dx
=e^{-\frac{\alpha^2}{2}}\Big[\frac{1}{\sqrt{2\pi}}
\int_\mathbb{R}e^{i\beta x}e^{-\frac{x^2}{2}}dx\Big]
=e^{-\frac{\alpha^2}{2}}e^{-\frac{\beta^2}{2}}
$$

\subsection{Vacuum characteristic functions of observables in $Heis(1,2)$}

In this section we use the $q$--projection method to give a derivation of 
the expression of the vacuum characteristic functions of observables 
in $\{1,p,q,q^2\}$ different from the one discussed in \cite{2}.
\begin{proposition}\label{prop1}
For all, $\alpha,\;\beta,\;\gamma\in\mathbb{R}$, one has
$$
\langle \Phi, e^{i\alpha(\sqrt{2}q)^2+i\beta(\sqrt{2}q)
+\gamma(\sqrt{2}q)}\Phi\rangle
=(1-2i\alpha)^{-\frac{1}{2}}e^{\frac{\gamma^2}{2(1-2i\alpha)}}
e^{-\frac{\beta^2}{2(1-2i\alpha)}}e^{i\frac{\beta\gamma}{1-2i\alpha}}
$$
\end{proposition}
\begin{proof}
Put
$$
\Psi_1(\beta):=\langle \Phi, 
e^{i\alpha(\sqrt{2}q)^2+i\beta(\sqrt{2}q)+\gamma(\sqrt{2}q)}\Phi\rangle
=\mathbb{E}(e^{i\alpha X^2+i\beta X+\gamma X})
$$
where $X$ is a normal gaussian random variable. Then, one has
\begin{eqnarray}\label{beta1}
\Psi_1'(\beta)&=&i\mathbb{E}(Xe^{i\alpha X^2+i\beta X+\gamma X})\\
&=&\frac{1}{\sqrt{2\pi}}\int_\mathbb{R}xe^{i\alpha x^2+i\beta x+\gamma x}
e^{-\frac{x^2}{2}}dx\nonumber
\end{eqnarray}
Taking the changes of variables
$$
u(x)=e^{i\alpha x^2+i\beta x+\gamma x},\;v'(x)=xe^{-\frac{x^2}{2}}
$$
Then, one gets
$$
u'(x)=(2i\alpha x+i\beta +\gamma)e^{i\alpha x^2+i\beta x+\gamma x}
\qquad;\quad v(x)=-e^{-\frac{x^2}{2}}
$$
This gives
\begin{eqnarray}\label{beta2}
\mathbb{E}(Xe^{i\alpha X^2+i\beta X+\gamma X})
=2i\alpha\mathbb{E}(Xe^{i\alpha X^2+i\beta X+\gamma X})
+(i\beta+\gamma)\mathbb{E}(e^{i\alpha X^2+i\beta X+\gamma X})
\end{eqnarray}
Therefore, from identities (\ref{beta1}) and (\ref{beta2}), one has
\begin{eqnarray*}
\Psi_1'(\beta)=\frac{i\gamma-\beta}{1-2i\alpha}\Psi_1(\beta)
\end{eqnarray*}
which yields
\begin{equation}\label{Beta}
\Psi_1(\beta)=C(\alpha,\gamma)e^{i\frac{\beta\gamma}{1-2i\alpha}}
e^{-\frac{\beta^2}{2(1-2i\alpha)}}
\end{equation}
where 
$$
C(\alpha,\gamma)=\Psi_1(0):=\Psi_2(\gamma)
=\mathbb{E}(e^{i\alpha X^2+\gamma X})
$$
Note that
\begin{eqnarray*}
\Psi_2'(\gamma)&=&\mathbb{E}(Xe^{i\alpha X^2+\gamma X})\\
&=&\frac{1}{\sqrt{2\pi}}\int_\mathbb{R}xe^{i\alpha x^2+\gamma x}
e^{-\frac{x^2}{2}}dx
\end{eqnarray*}
Taking the changes of variables
$$
h(x)=e^{i\alpha x^2+\gamma x},\;l'(x)=xe^{-\frac{x^2}{2}}
$$
It follows that
$$
h'(x)=(2i\alpha x+\gamma)e^{i\alpha x^2+\gamma x},\;l(x)=-e^{-\frac{x^2}{2}}
$$
Then, one obtains
\begin{eqnarray*}
\Psi_2'(\gamma)=2i\alpha\Psi_2'(\gamma)+\gamma\Psi_2(\gamma)
\end{eqnarray*}
This gives
\begin{equation}\label{Gamma}
\Psi_2(\gamma)=C(\alpha)e^{\frac{\gamma^2}{2(1-2i\alpha)}}
\end{equation}
where 
\begin{equation}\label{Alpha}
C(\alpha)=\Psi_2(0)=\mathbb{E}(e^{i\alpha X^2})=(1-2i\alpha)^{-\frac{1}{2}}
\end{equation}
Finally, using identities (\ref{Beta}), (\ref{Gamma}) and (\ref{Alpha}), 
one obtains
$$
\langle \Phi, e^{i\alpha(\sqrt{2}q)^2+i\beta(\sqrt{2}q)+\gamma(\sqrt{2}q)}
\Phi\rangle
=(1-2i\alpha)^{-\frac{1}{2}}e^{\frac{\gamma^2}{2(1-2i\alpha)}}
e^{-\frac{\beta^2}{2(1-2i\alpha)}}e^{i\frac{\beta\gamma}{1-2i\alpha}}
$$
\end{proof}
\begin{theorem}\label{characteristic}
For all $A,\,B,\,C\in\mathbb{R}$, one has
\begin{eqnarray*}
\langle \Phi,e^{it(A(\sqrt{2}q)^2+B(\sqrt{2}q)+C(\sqrt{2}p))}\Phi\rangle=(1-2itA)^{-\frac{1}{2}}e^{\frac{4C^2(A^2t^4+2iAt^3)-3|M|^2t^2}{6(1-2iAt)}}
\end{eqnarray*}
where $M=B+iC$.
\end{theorem}
\begin{proof}
We have
$$it(Aq^2+Bq+Cp)=itCp+itP'(q)$$
where $P(X)=\frac{1}{3}AX^3+\frac{1}{2}BX^2$. 
Then, Proposition \ref{prop2} implies that
\begin{equation}\label{decomp}
e^{it(Aq^2+Bq+Cp)}=e^{it\frac{P(q+tC)-P(q)}{tC}}e^{itCp}
\end{equation}
But, one has
\begin{equation}\label{decomp2}
\frac{P(q+tC)-P(q)}{tC}=Aq^2+(tAC+B)q+\frac{1}{2}tBC+\frac{1}{3}A(tC)^2
\end{equation}
Using (\ref{decomp}) and (\ref{decomp2}) for getting
\begin{eqnarray*}
e^{it(Aq^2+Bq+Cp)}
=e^{it\big(\frac{1}{3}A(tC)^2+\frac{1}{2}tBC\big)}
e^{it\big(Aq^2+(tAC+B)q\big)}e^{itCp}
\end{eqnarray*}
It follows that
\begin{eqnarray}\label{decomp3}
e^{it(A(\sqrt{2}q)^2+B(\sqrt{2}q)+C(\sqrt{2}p))}
=&&e^{it\big(\frac{4}{3}A(tC)^2+tBC\big)}
e^{it\big(A(\sqrt{2}q)^2+(2tAC+B)(\sqrt{2}q)\big)}\nonumber\\
&&e^{itC(\sqrt{2}p)}
\end{eqnarray}
Therefore from (\ref{decomp3}), (\ref{decomp4}) and (\ref{decomp4}), 
one has
\begin{eqnarray}\label{decomp5}
e^{it(A(\sqrt{2}q)^2+B(\sqrt{2}q)+C(\sqrt{2}p))}\Phi\!\!\!&=&\!\!\!
e^{it\big(\frac{4}{3}AC^2t^2+t(B+iC)C\big)}
e^{it\big(A(\sqrt{2}q)^2+(2ACt+(B+iC))(\sqrt{2}q)\big)}\Phi\nonumber\\
&=&e^{it\big(\frac{4}{3}AC^2t^2+tMC\big)}
e^{it\big(A(\sqrt{2}q)^2+(2ACt+M)(\sqrt{2}q)\big)}\Phi
\end{eqnarray}
where $M=B+iC$. Now, by taking
\begin{eqnarray*}
\alpha=At,\;\beta=2ACt^2+Bt,\;\gamma=-Ct
\end{eqnarray*}
and using Proposition \ref{prop1}, one gets
\begin{eqnarray}\label{sqrt}
\langle \Phi, e^{it\big(A(\sqrt{2}q)^2+(2ACt+M)(\sqrt{2}q)\big)}\Phi\rangle
=&&(1-2iAt)^{-\frac{1}{2}}e^{\frac{C^2t^2}{2(1-2iAt)}}
e^{-\frac{(2At^2+Bt)^2}{2(1-2iAt)}}\nonumber\\
&&e^{-i\frac{Ct(2ACt^2+Bt)}{1-2iAt}}
\end{eqnarray}
Finally, identities (\ref{decomp5}) and (\ref{sqrt}) imply that
$$
\langle \Phi,e^{it(A(\sqrt{2}q)^2+B(\sqrt{2}q)+C(\sqrt{2}p))}\Phi\rangle
=(1-2itA)^{-\frac{1}{2}}e^{\frac{4C^2(A^2t^4+2iAt^3)-3|M|^2t^2}{6(1-2iAt)}}
$$
This ends the proof.
\end{proof}\\
Using together Proposition \ref{N=2} and Theorem 
\ref{characteristic}, we prove the following theorem.
\begin{theorem}
For all $\alpha_i,\beta_i,\gamma_i\in\mathbb{R},\;i=1,2,$ we have
\begin{eqnarray*}
&&\langle e^{it(\alpha_1(\sqrt{2}q)^2+\beta_1(\sqrt{2}q)
+\gamma_1(\sqrt{2}p))}\Phi,e^{it(\alpha_2(\sqrt{2}q)^2
+\beta_2(\sqrt{2}q)+\gamma_2(\sqrt{2}p))}\Phi\rangle\\
&&=(1-2i(\alpha_2-\alpha_1)t)^{-\frac{1}{2}}
e^{\frac{Lt^4+Z_1t^3+Z_2t^2}{6(1-2i(\alpha_2-\alpha_1)t)}}
\end{eqnarray*}
where
\begin{eqnarray*}
L&=&-4(\alpha_1\gamma_2-\alpha_2\gamma_1)
\big[3(\alpha_1\gamma_2-\alpha_2\gamma_1)+2(\alpha_2-\alpha_1)(\gamma_1+\gamma_2)\big] \\
&&+4(\gamma_2-\gamma_1)^2(\alpha_2-\alpha_1)^2\\
Z_1&=&12\big[(\alpha_2-\alpha_1)(\beta_1\gamma_2-\gamma_1\beta_2)-(\beta_2-\beta_1)(\alpha_1\gamma_2-\alpha_2\gamma_1)\big]\\
&&+4i\big[2(\alpha_2-\alpha_1)(\gamma_2-\gamma_1)^2+(\gamma_1+\gamma_2)(\alpha_2\gamma_1-\alpha_1\gamma_2)\big]\\
Z_2&=& -3(\gamma_2-\gamma_1)^2+6i(\beta_1\gamma_2-\gamma_1\beta_2)  
\end{eqnarray*}
\end{theorem}
\begin{proof}
We have
\begin{eqnarray}\label{eqnn}
&&\langle e^{it(\alpha_1(\sqrt{2}q)^2+\beta_1(\sqrt{2}q)
+\gamma_1(\sqrt{2}p))}\Phi,e^{it(\alpha_2(\sqrt{2}q)^2
+\beta_2(\sqrt{2}q)+\gamma_2(\sqrt{2}p))}\Phi\rangle\nonumber\\
&&=\langle \Phi,e^{-it(\alpha_1(\sqrt{2}q)^2+\beta_1(\sqrt{2}q)
+\gamma_1(\sqrt{2}p))}e^{it(\alpha_2(\sqrt{2}q)^2+\beta_2(\sqrt{2}q)
+\gamma_2(\sqrt{2}p))}\Phi\rangle
\end{eqnarray}
Put
\begin{eqnarray*}
P'_1(x)=-t\sqrt{2}\beta_1x-2t\alpha_1x^2\\
P'_2(x)=t\sqrt{2}\beta_2x+2t\alpha_2x^2
\end{eqnarray*}
Then, Proposition \ref{N=2} implies that
\begin{equation}\label{ABC}
(-t\sqrt{2}\gamma_1,P'_1)(t\sqrt{2}\gamma_2,P_2')
=(t\sqrt{2}(\gamma_2-\gamma_1),Q')
\end{equation}
where $Q'(x)=C+Bx+Ax^2$ with
\begin{eqnarray*}
C&=&(\beta_1\gamma_2-\gamma_1\beta_2)t^2+\frac{2}{3}
(\gamma_1+\gamma_2)(\alpha_2\gamma_1-\alpha_1\gamma_2)t^3\\
B&=&\sqrt{2}(\beta_2-\beta_1)t
+2\sqrt{2}(\alpha_1\gamma_2-\alpha_2\gamma_1)t^2\\
A&=&2(\alpha_2-\alpha_1)t
\end{eqnarray*}
Therefore, from identities (\ref{eqnn}) and (\ref{ABC}), one obtains
\begin{eqnarray*}
&&\langle e^{it(\alpha_1(\sqrt{2}q)^2+\beta_1(\sqrt{2}q)+\gamma_1(\sqrt{2}p))}\Phi,e^{it(\alpha_2(\sqrt{2}q)^2+\beta_2(\sqrt{2}q)+\gamma_2(\sqrt{2}p))}\Phi\rangle\\
&&=e^{i\big((\beta_1\gamma_2-\gamma_1\beta_2)t^2+\frac{2}{3}(\gamma_1+\gamma_2)(\alpha_2\gamma_1-\alpha_1\gamma_2)t^3\big)}\\
&&\;\;\;\;\;\;\;\langle\Phi,e^{it\Big((\alpha_2-\alpha_1)(\sqrt{2}q)^2+\big((\beta_2-\beta_1)+2(\alpha_1\gamma_2-\alpha_2\gamma_1)t\big)(\sqrt{2}q)+(\gamma_2-\gamma_1)(\sqrt{2}p)\Big)}\Phi\rangle
\end{eqnarray*}
Now, from Theorem \ref{characteristic}, one gets
\begin{eqnarray}\label{Mm}
&&\langle e^{it(\alpha_1(\sqrt{2}q)^2+\beta_1(\sqrt{2}q)+\gamma_1(\sqrt{2}p))}\Phi,e^{it(\alpha_2(\sqrt{2}q)^2+\beta_2(\sqrt{2}q)+\gamma_2(\sqrt{2}p))}\Phi\rangle\nonumber\\
&&=e^{i\big((\beta_1\gamma_2-\gamma_1\beta_2)t^2+\frac{2}{3}(\gamma_1+\gamma_2)(\alpha_2\gamma_1-\alpha_1\gamma_2)t^3\big)}(1-2i(\alpha_2-\alpha_1)t)^{-\frac{1}{2}}\nonumber\\
&&\;\;\;\,\;e^{\frac{4(\gamma_2-\gamma_1)^2\big((\alpha_2-\alpha_1)^2t^4+2i(\alpha_2-\alpha_1)t^3\big)-3|M|^2t^2}{6(1-2i(\alpha_2-\alpha_1)t)}}
\end{eqnarray}
where 
$$M=(\beta_2-\beta_1)+2(\alpha_1\gamma_2-\alpha_2\gamma_1)t+i(\gamma_2-\gamma_1)$$
Finally, using the identity (\ref{Mm}), it is easy to show that
\begin{eqnarray*}
&&\langle e^{it(\alpha_1(\sqrt{2}q)^2+\beta_1(\sqrt{2}q)+\gamma_1(\sqrt{2}p))}\Phi,e^{it(\alpha_2(\sqrt{2}q)^2+\beta_2(\sqrt{2}q)+\gamma_2(\sqrt{2}p))}\Phi\rangle\\
&&=(1-2i(\alpha_2-\alpha_1)t)^{-\frac{1}{2}}e^{\frac{Lt^4+Z_1t^3+Z_2t^2}{6(1-2i(\alpha_2-\alpha_1)t)}}
\end{eqnarray*}
where $L, Z_1,Z_2$ are given above.
\end{proof}

\subsection{Vacuum moments of observables in $\{1,p,q,q^2\}$}\label{3.3}

In this section we use the results of the preceeding section to deduce
the expression of the vacuum moments of observables in $\{1,p,q,q^2\}$.
The form of these moments was not known and may be used to throw some light 
on the still open problem of finding the explicit expression of 
the probability distributions of these observables.
\begin{theorem}
Define, for $A,B,C\in\mathbb R$
$$
X:=A(\sqrt{2}q)^2+B(\sqrt{2}q)+C(\sqrt{2}p)
$$
Then
$$
\langle \Phi,X^n\Phi\rangle
=\sum_{i_1+2i_2+\dots+ki_k=n}\frac{2^{3n}n!}{i_1!\dots i_k!}w_1^{i_1}\dots w_k^{i_k}
$$
where
$$
w_k=\frac{A^k}{2k}-\frac{A^{k-2}}{4}\gamma\chi_{\{k\geq2\}}-\frac{A^{k-3}}{8}\beta\chi_{\{k\geq3\}}+\frac{A^{k-4}}{16}\alpha\chi_{\{k\geq4\}}
$$
with
$$
\alpha=\frac{2}{3}A^2C^2,\;\beta=\frac{4}{3}AC^2,\;
\gamma=-\frac{1}{2}|M|^2=-\frac{1}{2}(B^2 +C^2)
$$
\end{theorem}
\begin{proof}
Recall that
$$\mathbb{E}(e^{itX})=e^{\varphi(t)}$$
where 
$$\varphi(t)=-\frac{1}{2}\ln(1-2iAt)+\frac{\alpha t^4+i\beta t^3+\gamma t^2}{1-2iAt}$$
with $\alpha,\beta$ and $\gamma$ are given above. Hence, one has
\begin{equation}\label{n-moment}
\mathbb{E}(X^n)=\frac{1}{i^n}\Big(\frac{d}{dt}\Big)^n\Big|_{t=0}e^{\varphi(t)}
\end{equation}
Now we introduce the following formula (cf \cite{3})
\begin{eqnarray}\label{bourbaki}
\frac{d^n}{dt^n}e^{\varphi(t)}
=\sum_{i_1+2i_2+\dots+ki_k=n}\frac{2^{2n}n!}{i_1!\dots 
i_k!}\Big(\frac{\varphi^{(1)}(t)}{1!}\Big)^{i_1}\dots
\Big(\frac{\varphi^{(k)}(t)}{k!}\Big)^{i_k}e^{\varphi(t)}
\end{eqnarray}
Put 
\begin{eqnarray*}
&&\varphi_1(t)=\alpha t^4+i\beta t^3+\gamma t^2,\;\varphi_2(t)=(1-2iAt)^{-1},\\
&&\varphi_3(t)=-\frac{1}{2}\ln(1-2iAt),\;g(t)=\varphi_1(t)\varphi_2(t)
\end{eqnarray*}
Note that
\begin{eqnarray}\label{derivative}
\varphi_2^{(k)}(t)=(2iA)^kk!(1-2iAt)^{-k-1},\varphi_3^{(k)}(t)
=\frac{1}{2}(2iA)^k(k-1)!(1-2iAt)^{-k}
\end{eqnarray}
Then, one gets
\begin{eqnarray*}
g^{(k)}(t)=\sum_{h=0}^k\pmatrix{k\cr h\cr}\varphi_1^{(h)}(t)
\varphi_2^{(k-h)}(t)
\end{eqnarray*}
Because $\varphi_1^{(h)}(t)=0$ for all $h\geq5$ and $\varphi_1(0)=\varphi_1'(0)=0$, one gets
\begin{eqnarray}\label{convention}
g^{(k)}(0)=\frac{k!}{(k-2)!}\gamma\varphi_2^{(k-2)}(0)
+i\frac{k!}{(k-3)!}\beta\varphi_2^{(k-3)}(0)+\frac{k!}{(k-4)!}\alpha
\varphi_2^{(k-4)}(0)
\end{eqnarray}
where by convention $\varphi_2^{(i-j)}(0)=0$ if $i<j$. 
Therefore, identities (\ref{derivative}) and (\ref{convention}) imply that
\begin{eqnarray*}
g^{(k)}(0)=k!(2iA)^{k-2}\gamma\chi_{\{k\geq2\}}+ik!(2iA)^{k-3}\beta\chi_{\{k\geq3\}}+k!(2iA)^{k-4}\alpha\chi_{\{k\geq4\}}
\end{eqnarray*}
Then, for all $k\geq1$, one obtains
\begin{eqnarray*}
\frac{\varphi^{(k)}(0)}{k!}&=&\varphi_3^{(k)}(0)+g^{(k)}(0)\\
&=&(2i)^k\Big(\frac{A^k}{2k}-\frac{A^{k-2}}{4}\gamma\chi_{\{k\geq2\}}
-\frac{A^{k-3}}{8}\beta\chi_{\{k\geq3\}}
+\frac{A^{k-4}}{16}\alpha\chi_{\{k\geq4\}}\Big)\\
&=&(2i)^kw_k
\end{eqnarray*}
where
$$
w_k=\frac{A^k}{2k}-\frac{A^{k-2}}{4}\gamma\chi_{\{k\geq2\}}
-\frac{A^{k-3}}{8}\beta\chi_{\{k\geq3\}}
+\frac{A^{k-4}}{16}\alpha\chi_{\{k\geq4\}}
$$
Thus, from identity (\ref{bourbaki}), one has 
\begin{eqnarray*}
\Big(\frac{d}{dt}\Big)^n\Big|_{t=0}e^{\varphi(t)}&=&\sum_{i_1+2i_2+\dots+ki_k=n}\frac{2^{2n}n!}{i_1!\dots i_k!}(2i)^{i_1+\dots+ki_k}w_1^{i_1}\dots w_k^{i_k}\\
&=&i^n\sum_{i_1+2i_2+\dots+ki_k=n}\frac{2^{2n}2^nn!}{i_1!\dots i_k!}w_1^{i_1}\dots w_k^{i_k}
\end{eqnarray*}
Finally, by using identity (\ref{n-moment}) the result of the above theorem holds true.
\end{proof}

\subsection{Vacuum characteristic functions of observables in $Heis(1,n)$}

From (\ref{*2}) and  (\ref{decomp4}) we deduce that
\begin{eqnarray*}
\langle \Phi , e^{iwp+iuP'(q)}\Phi \rangle &=&\langle \Phi ,e^{iu\frac{P(q+w)-P(q)}{w}}e^{iwp}\Phi \rangle \\
&=&\langle \Phi ,e^{iu\frac{P(q+w)-P(q)}{w}}e^{iwp}\Phi \rangle \\
&=&e^{-w^2/2}\langle \Phi ,e^{iu\frac{P(q+w)-P(q)}{w}}e^{-wq}\Phi \rangle 
\end{eqnarray*}
Thus the q-projection method reduces the problem to an integral of the form
\begin{eqnarray*}
\langle \Phi , e^{iQ(q)}\Phi \rangle 
\end{eqnarray*}
where $Q$ is a polynomial.\\

{\bf Acknowledgements.}\\
LA acknowledges support by the RSF grant 14-11-00687.

\end{document}